\documentclass{amsart}
\usepackage[dvips]{graphics}
\usepackage{amsthm,amssymb}

\usepackage[dvips]{color}

\theoremstyle{definition}

\theoremstyle{remark}

\numberwithin{equation}{section}

\begin{document}

\title{The Dehn filling space of a certain hyperbolic 3-orbifold}

\author{Sadayoshi Kojima}
\address{Department of Mathematical and Computing Sciences, 
	Tokyo Institute of Technology, 
	Meguro, Tokyo 152-8552 Japan}
\email{sadayosi@is.titech.ac.jp}
\thanks{The first author was supported in part by JSPS Grant 12440015.}

\author{Shigeru Mizushima}
\address{Department of Mathematical and Computing Sciences, 
	Tokyo Institute of Technology, 
	Meguro, Tokyo 152-8552 Japan}
\email{mizusima@is.titech.ac.jp}




\begin{abstract}
We construct the first example of a ``one-cusped'' hyperbolic 3-orbifold 
for which we see the true shape of the space of 
hyperbolic Dehn fillings. 
\end{abstract}

\maketitle


%
%

\section{Introduction}

The hyperbolic Dehn filling theory was established by 
W. Thurston in his original note  \cite{Thu}.  
Besides the general theory, 
he extensively analyzed the space of Dehn fillings of 
the figure eight knot complement in terms of the decomposition by 
two ideal tetrahedra. 
The analysis leaded us to 
extremely fascinated research activities on Dehn fillings 
in the last 25 years.  

One of unsolved problems in his analysis was to determine 
the global shape of the Dehn filling space. 
Since then, 
there have been many deep researches such as 
\cite{Hodgson, Dowty} discussing this problem in fact, 
however it is still mysterious.  
Even, 
though W. Neumann and A. Reid  \cite{NR}  have determined the 
true shape 
for the Whitehead link complement with one cusp component being complete, 
there seems to be no concrete examples of ``one-cusped'' 
hyperbolic manifolds 
or orbifolds for which we see the true shape 
of the Dehn filling space.   
In this paper, 
we construct hopefully the first such example 
of an orbifold, 
regretfully rather than a manifold. 

The example we present here is 
topologically  ${\mathbb Z}_2 \times {\mathbb Z}_2$-covered 
as an orbifold by 
the complement of a seven component link in the 
connected sum of three copies of  ${\bf S}^2 \times {\bf S}^1$.  
Based on the study of one-circle packings on complex affine tori 
by the second author in \cite{Miz},  
we construct the example, 
which we will denote by  $N$,  
and its all possible Dehn filling deformations 
in the following three sections. 
Then using careful analysis of deformations of such tori 
in  \cite{Miz} again,  
we determine the global shape of the Dehn filling space of  $N$  
in subsequent two sections.

%
%

\section{Construction}

\begin{figure}[hbt]
  \begin{center}
    \begin{picture}(100,100)
       \put(0,0){\scalebox{0.5}{\includegraphics{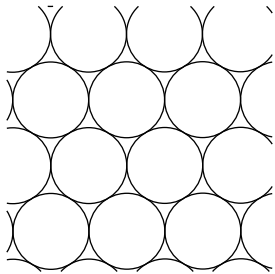}}}
    \end{picture}
    \caption[]{Hexagonal packing}
    \label{Fig:HexagonalPacking}
  \end{center}
\end{figure}

We start with the {\it hexagonal packing}  $\mathcal{P}$  
on the complex plane  $\mathbb{C}$  by equi-radii circles, 
see Figure \ref{Fig:HexagonalPacking}.  
The set of euclidean translations which leave  $\mathcal{P}$  
invariant forms a group  $\Gamma$  acting freely on  $\mathbb{C}$. 
It is isomorphic to  $\mathbb{Z} \times \mathbb{Z}$  and 
the quotient is a particular euclidean torus admitting 
a cyclic symmetry of order  $6$,  
which we call the {\it hexagonal torus}. 
$\mathcal{P}$  then descends to what we call a one-circle 
packing on the hexagonal torus. 

Let us regard the complex plane  $\mathbb{C}$  as 
the boundary of the upper half 
space model  $\mathbb{H}^3$  of the hyperbolic 3-space. 
$\mathcal{P}$  together with the dual packing  $\mathcal{P}^*$  which 
consists of circumscribed circles of interstices in  $\mathcal{P}$  
fills  $\mathbb{C}$.  
Each member of both  $\mathcal{P}$  and  $\mathcal{P}^*$  bounds 
a hemisphere in  $\mathbb{H}^3$  and 
such a hemisphere bounds a hemiball facing  $\mathbb{C}$.  
Cutting off those hemiballs from  $\mathbb{H}^3$, 
we obtain a region  $L$  with ideal polygonal boundary  in  $\mathbb{H}^3$. 
The intersections of  $\mathcal{P}$  and $\mathcal{P}^*$  correspond
to ideal vertices of  $L$  with rectangular section. 

The group  $\Gamma$  acts properly discontinuously on 
$\mathbb{H}^3$  by the Poincar\'e extensions and 
the region  $L$  is invariant under the action.  
Hence taking quotient, 
we obtain a hyperbolic 3-manifold  $N$  with 
ideal polygonal boundary and a cusp. 
The polygonal boundary consists of 
two ideal triangles and one ideal hexagon, 
where they intersect orthogonally. 

Take the double of  $N$  along triangular faces, 
and then take again the double of the result along 
the boundary which now consists of two ideal hexagons. 
This double doubling construction gives us  
a hyperbolic 3-manifold  $M$  with seven cusps among which 
4 are old and 3 new. 
It is obvious by the construction that  $M$  admits 
a  $\mathbb{Z}_2 \times \mathbb{Z}_2$-symmetry generated 
by two reflections associated to the doubling.  
$N$  is the quotient of  $M$  by this symmetry 
and therefore we regard  $N$  as an orbifold from now on, 
and study its Dehn filling deformations as an orbifold.  

\begin{figure}[hbt]
  \begin{center}
    \begin{picture}(141,123)
       \put(0,0){\scalebox{0.5}{\includegraphics{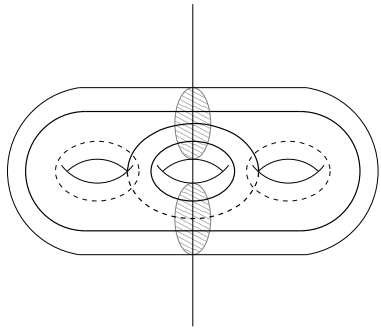}}}
    \end{picture}
    \caption[]{First double}
    \label{Fig:FirstDouble}
  \end{center}
\end{figure}

Topologically  $N$  is a solid torus with a core loop and 
three points on the boundary removed.  
To see the effect on the first double, 
stretch the removed vertices on  $\partial N$  to arcs appropriately 
so that they become circles after doubling. 
Then by first doubling, 
we obtain a genus 3 handle body 
with two core loops and three loops on the boundary removed. 
This is depicted in Figure \ref{Fig:FirstDouble}, 
where the shaded disks represent triangular faces along which 
we took the double. 

\begin{figure}[hbt]
  \begin{center}
    \begin{picture}(145,110)
       \put(0,0){\scalebox{0.5}{\includegraphics{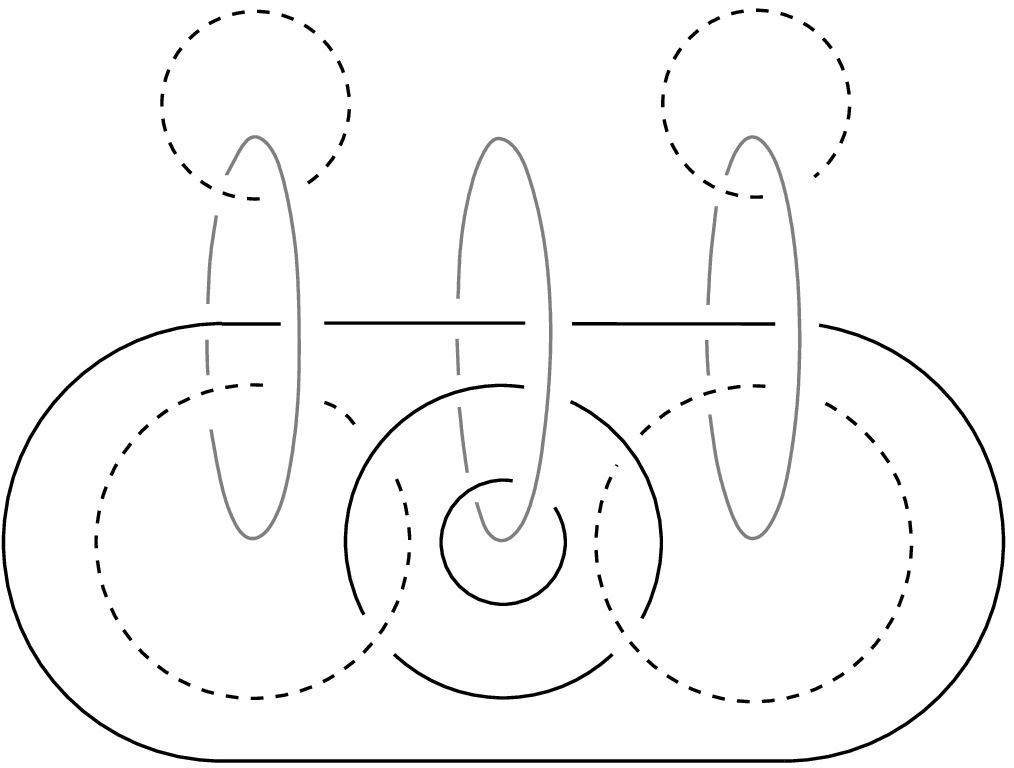}}}
       \put(45,72){$0$}
       \put(81,72){$0$}
       \put(117,72){$0$}
    \end{picture}
    \caption[]{Framed link representation of  $M$}
    \label{Fig:FramedLink}
  \end{center}
\end{figure}

A framed link representation of  $M$, 
the double of Figure \ref{Fig:FirstDouble}, 
can be immediately described by Figure \ref{Fig:FramedLink}. 

%
%

\section{Circle Packings on Complex Affine Tori}

We will construct Dehn filling deformations of  $N$  using 
deformations of the hexagonal packing on  $\mathbb{C}$  
with the same combinatorial type.  
The hexagonal packing descends to a one-circle packing 
on the hexagonal torus. 
Conversely, 
a one-circle packing on the torus combinatorially equivalent 
to one on the hexagonal torus is universal covered by 
a packing on the plane equivalent to the hexagonal packing. 
However since the hexagonal packing is rigid in euclidean geometry, 
there are no other one-circle packings on euclidean tori 
than one on the hexagonal torus.  
Hence to get deformations, 
we will work with one-circle packings not on  
euclidean but rather on complex affine tori.  
This section is to review what has been known about complex affine 
tori and circle packings on them 
according to  \cite{Gun, Miz}. 

A {\it complex affine structure} on the torus is by definition 
a collection of local charts on  $\mathbb{C}$  whose 
transition function is 
a restriction of a complex affine transformation, 
so in particular it defines a complex structure.  
The uniformization assigning the conformal class 
to each affine torus defines a map of 
the space  $\mathcal{A}$  of marked complex affine tori  
to the Teichm\"uller space  $\mathcal{H}$  of the torus.  
It is a natural complex line bundle projection, 
where the fiber is parameterized by explicit description of   
developing maps of affine structures under some normalization 
as we will see below.  
$\mathcal{H}$  is identified with the upper half plane  $\mathbb{H}^2$  
in  $\mathbb{C}$  and hence  $\mathcal{A}$  is homeomorphic to 
the product  $\mathbb{H}^2 \times \mathbb{C}$. 
We call the first factor a {\it Teichm\"uller parameter}  
$\omega \in \mathbb{H}^2$  and 
the second an {\it affine parameter} $c \in \mathbb{C}$.    
Complex affine tori sharing the same $\omega$ are conformally equivalent.
We denote by  $T_{\omega,c}$  the marked complex affine torus 
associated to  $(\omega, c) \in \mathbb{H}^2 \times \mathbb{C}$.  
Also we denote by  $\alpha, \, \beta$  specified 
simple closed curves generating 
the fundamental group of the torus for marking.   

The zero section of 
the line bundle $\mathbb{H}^2 \times \mathbb{C} \to \mathbb{H}^2$  
precisely consists of euclidean structures.  
In this case, 
the developing map 
\begin{equation*} 
	{\rm dev}:\widetilde{T}_{\omega,0} \rightarrow \mathbb{C}, 
\end{equation*}
is a homeomorphism, 
and the images of  $\alpha$ and $\beta$  by the holonomy $\rho$  
are translations defined for each  $z \in \mathbb{C}$  by 
\begin{equation*} 
	\begin{cases} 
		\quad \rho(\alpha):z\mapsto z+1 \\ 
		\quad \rho(\beta):z\mapsto z+\omega 
	\end{cases} 
\end{equation*}  
after appropriate normalization.

\begin{figure}[hbt]
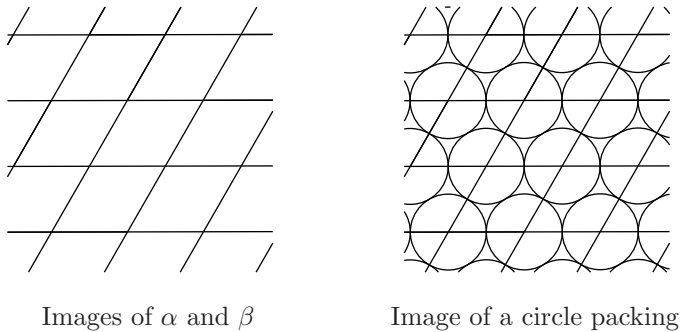

  \begin{center}
    \begin{picture}(250,120)
       \put(0,20){\scalebox{0.5}{\includegraphics{Fundamental-Region-Euclidean.eps}}}
	\put(13,0){Images of $\alpha$ and $\beta$}
       \put(150,20){\scalebox{0.5}{\includegraphics{Packing-Fundamental-Region-Euclidean2.eps}}}
	\put(145,0){Image of a circle packing}
    \end{picture}
    \caption[]{Developed image of torus with euclidean structure}
    \label{Fig:EuclideanStr}
  \end{center}
\end{figure}

In affine but non-euclidean case,
identifying the universal cover  $\widetilde{T}_{\omega,c}$  with 
$\widetilde{T}_{\omega,0} \cong \mathbb{C}$  by the uniformization, 
the developing map is described by an exponential map 
\begin{equation*} 
	{\rm dev}: z \mapsto e^{cz} 
\end{equation*} 
under the identification of  $\mathcal{A}$  with  
$\mathbb{H}^2 \times \mathbb{C}$. 
Here $\omega$  and  $c$  are the Teichm\"uller  
and affine parameters of the complex affine torus 
in question.  
The images of  $\alpha, \beta$  by the holonomy $\rho$  are  
obviously similarities 
\begin{equation*} 
	\begin{cases} 
		\quad \rho(\alpha):z \mapsto e^{c}z \\ 
		\quad \rho(\beta):z \mapsto e^{c\omega}z 
	\end{cases} 
\end{equation*} 
which fix the origin.  

\begin{figure}[hbt]
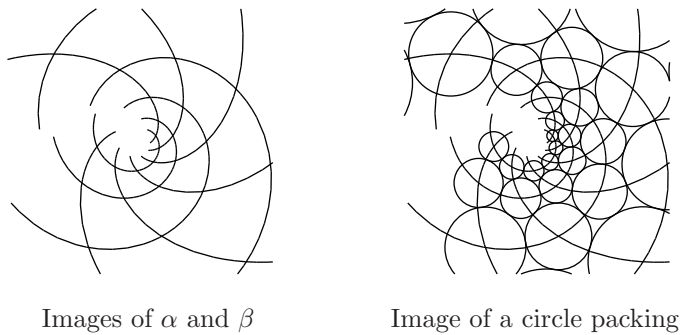

  \begin{center}
    \begin{picture}(250,120)
       \put(0,20){\scalebox{0.5}{\includegraphics{Affine-Curves.eps}}}
	\put(13,0){Images of $\alpha$ and $\beta$}
       \put(150,20){\scalebox{0.5}{\includegraphics{Affine-CirclesCurves.eps}}}
	\put(145,0){Image of a circle packing}
    \end{picture}
    \caption[]{Developed image of torus with complex affine structure}
    \label{Fig:AffineStr}
  \end{center}
\end{figure}

Since a complex affine transformation is a similarity,  
it sends a circle on  $\mathbb{C}$  to a circle.  
Hence we can still discuss about circle packings on complex affine tori. 
The second author showed in  \cite{Miz} that 
the set of complex affine tori which admit 
a one-circle packing  forms 
a real two-dimensional moduli space  $\mathcal{M}$  
in  $\mathbb{H}^2 \times \mathbb{C}$. 
In constructing the moduli space, 
the marking for the packing was specified so that 
the edges of the nerve of 
the packing are homotopic to  $\alpha$, $\beta$  
and  $\beta\alpha^{-1}$.
The moduli space  $\mathcal{M}$  contains only one euclidean torus 
associated to  $\omega_0 = e^{\pi i/3}$,  
which is the hexagonal one. 
Moreover, 
the restriction of 
the uniformization $\mathbb{H}^2 \times \mathbb{C} \to \mathbb{H}^2$  
to  $\mathcal{M}$  is onto and 
two to one except at the hexagonal torus. 
The affine parameters of a generic fiber consist of nonzero  $c$  and 
its negative  $-c$.  

%
%

\section{Deformations} 

To each pair  $(\omega,c)$  on the moduli space  $\mathcal{M}^\prime=\mathcal{M}-\{(\omega_0,0)\}$, 
we construct a Dehn filling deformation of  $N$  as an orbifold 
based on the construction of  $N$. 
To see this, 
we start with a complex affine torus  $T_{\omega, c}$  packed by one circle  
combinatorially equivalent to the one-circle packing 
on  $T_{\omega_0, 0}$.  
The developing map of  $T_{\omega, c}$  is an exponential map 
and $\mathbb{Z}$-covers  $\mathbb{C} - \{0\}$.  
Thus if we let  $X$  be the universal cover of  $\mathbb{C}-\{ 0 \}$,  
it canonically lifts to a homeomorphism  $:\widetilde{T}_{\omega,c} \to X$.  
By mapping a circle packing on  $\widetilde{T}_{\omega,c}$  to  $X$, 
we obtain a circle packing  $\mathcal{P}_{\omega,c}$  
on  $X$  combinatorially equivalent to the hexagonal packing 
on  $\mathbb{C}$.   
Also the holonomy lifts uniquely to an injective homomorphism of 
$\pi_1(T_{\omega,c}) \cong \mathbb{Z} \times \mathbb{Z}$  to 
the group of transformations of  $X$.  
We denote its image by  $\Gamma_{\omega,c}$.  

Then we let  $\ell$  denote the geodesic in  $\mathbb{H}^3$  
joining  $0$  and  $\infty$. 
The universal cover  $Y$  of  $\mathbb{H}^3 - \ell$  will be 
a substitution of  $\mathbb{H}^3$  in the euclidean case.  
In fact, 
cutting off hemiballs in  $Y$  bounded by the members of  
a deformed hexagonal packing  $\mathcal{P}_{\omega,c}$  and 
its dual  $\mathcal{P}_{\omega,c}^*$  on  $X$, 
we obtain a region  $L_{\omega,c}$  with ideal polygonal boundary in  $Y$.  
Also since the action of  $\Gamma_{\omega,c}$  on  $X$  extends to 
a properly discontinuous action on  $Y$  which leaves 
$L_{\omega,c}$  invariant, 
by taking a quotient we obtain a hyperbolic 3-manifold 
with ideal polygonal boundary 
and incomplete end near  $\ell$.  
We denote its completion by  $N_{\omega,c}$.  
It carries a Dehn surgery type singularity by definition, 
see  \cite{Thu}, 
and we thus obtain to each  $(\omega,c)$  in  $\mathcal{M}^\prime$  
a Dehn filling deformation of  $N = N_{\omega_0,0}$.  

On the other hand, 
let  $N_*$  be a Dehn filled deformation of  $N$  as an orbifold and  
$R$  a manifold obtained from  $N_*$  by removing 
the points attached by the completion. 
$R$  is homeomorphic to the torus times  $[0,1)$  with 
three points on the boundary removed.  
The boundary of  $R$  consists of two ideal triangles and 
one ideal hexagon, 
where they intersect orthogonally. 
Choose a developing map  $D$  of  $R$  so that a neighborhood of 
an incomplete end winds  $\ell$  in   $\mathbb{H}^3$.  
By this normalization, 
the holonomy of  $R$  maps  
$\pi_1(R) \cong \mathbb{Z} \times \mathbb{Z}$  to 
transformations fixing  $0$  and  $\infty$,  
the image of the holonomy contains an element of infinite order 
by the definition of the Dehn filling deformation. 

$R$  is not complete, 
but not so wild. 
Let  $U$  be an equidistant neighborhood of an incomplete end 
and  $B$  its complement. 
It is obvious that any point in  $R$  admits 
the unique shortest path to  $\partial U = \partial B$.   
Thus  $U$  is developed to an equidistant neighborhood of  $\ell$  
and  $B$  its complement.  
In particular, 
the image of the developing map  $D$  misses  $\ell$, 
and  $D$  is in fact a map of the universal cover  $\widetilde{R}$  
into  $\mathbb{H}^3 - \ell$. 
Then  $D$  lifts to a map  $\widetilde{D}$  of 
the universal cover  $\widetilde{R}$  to  $Y$, 
the universal cover of  $\mathbb{H}^3 - \ell$.  
Accordingly, 
the holonomy of  $R$  lifts to a homomorphism to the group of 
transformations of  $Y$.  
We denote its image by  $\Gamma_*$.  

Let  $Z$  be the preimage of the boundary of the equidistant neighborhood 
of  $\ell$  in  $Y$.  
$\Gamma_*$  acts properly discontinuously and freely on  $Z$  by 
the definition of the Dehn filling deformation,  
and  $Z$  bijectively corresponds to the preimage of  $Z$  
by  $\widetilde{D}$.  
This shows that  $\widetilde{D}$  is injective at least 
on a neighborhood of the preimage of  $Z$  by  $\widetilde{D}$.  
We will show that  $\widetilde{D}$  is globally injective. 
Suppose we have two points in  $\widetilde{R}$  mapped 
by  $\widetilde{D}$  to a point  $y \in Y$.  
$y$  admits the unique shortest path to the point  $y_0$  on  $Z$.  
Since the preimage of  $Z$  by  $\widetilde{D}$  bijectively 
corresponds to the preimage of  $\partial U = \partial B$  
in  $\widetilde{R}$, 
the preimage of  $y_0$  consists of a single point  $n_0$.  
Hence two points in question lie on the orthogonal path to 
the preimage of  $\partial U = \partial B$  through  $n_0$  on 
the same side with the same distance to  $n_0$.  
Thus they must be the same.  

Then the ideal hexagonal faces on the boundary of  $\widetilde{R}$  
in the developed image define a packing on  $X$, 
the universal cover of  $\mathbb{C} - \{0\}$,  combinatorially 
equivalent to the hexagonal one.  
Moreover the packing is invariant under the induced action of 
$\Gamma_*$  on  $X$.
Hence it descends to a one-circle packing on some affine torus 
$T_{\omega,c}$  parameterized by some pair  $(\omega, c) \in \mathcal{M}^\prime$, 
and  $N_*$  is isometric to  $N_{\omega,c}$.

%
%

\section{Dehn Filling Coefficients}

Following \cite{Miz}, 
we review how to find a pair  $(\omega,c)$  lying on 
the moduli space  $\mathcal{M}^\prime$.  
If  $(\omega,c)$  lies on  $\mathcal{M}^\prime$, 
a simple observation shows that  $\omega$  and  $c$  must satisfy 
the identities, 
\begin{eqnarray}
 \frac{\cos {\rm Im} (c/2)}{\cosh {\rm Re} (c/2)} =
 \frac{\cos {\rm Im} (c \omega / 2)}{\cosh {\rm Re} (c \omega/2)} =
 \frac{\cos {\rm Im} (c ( \omega - 1 ) /2)}
		{\cosh {\rm Re} (c ( \omega - 1 ) /2)}.
\label{Eq:ConditionForAdmittingPacking}
\end{eqnarray}
Moreover, the inequalities 
\begin{equation*} 
	|\mathrm{Im}(c)| < \pi, \quad 
	|\mathrm{Im}(c \omega)| < \pi, \quad  
	|\mathrm{Im}( c ( \omega - 1)) | < \pi
\end{equation*} 
must be satisfied by obvious reason.  
Let $f : \mathcal{C} \to [0,1)$  be a smooth function defined by
\begin{equation*} 
	f(z)=1-\frac{\cos {\rm Im} (z/2)}{\cosh {\rm Re}(z/2)}, 
\end{equation*} 
where 
\begin{equation*} 
	\mathcal{C} = \{z\in\mathbb{C}~|~|{\rm Im}(z)|<\pi\}
\end{equation*} 
is the largest possible region for  $c$  to vary in  $\mathcal{M}$.
The identity (\ref{Eq:ConditionForAdmittingPacking}) shows 
that $c$, $c\omega$ and $c(\omega-1)$ lie on 
the same level set  $L_s=\{z\in \mathcal{C}~|~f(z)=s\}$.
When  $s = 0$,  
$L_0$  consists of a single point 
corresponding to the hexagonal torus by simple computation. 
For each  $0<s<1$,  
$L_s$  is topologically a simple loop around the origin,
has symmetries about the real and imaginary axis, 
and bounds a strictly convex interior.

For $c \in \mathcal{C}^\prime=\mathcal{C}-\{0\}$,
the origin  $0$, $c$, 
$c\omega$ and  $c(\omega-1)$  form a parallelogram 
by obvious reason in plane geometry.
Conversely, 
if $c \in \mathcal{C}^\prime $ is given and let  $L_s$  be the 
level set on which  $c$  lies,  
then there are unique two points  $c_1, c_2$  on  $L_s$  
so that  $c, \, c_1, c_2$   lie in counterclockwise order on  $L_s$,  
and they together with the origin form a parallelogram. 
Then we obtain a pair  $(\omega,c)$  in the moduli space  $\mathcal{M}^\prime$  by 
letting  $\omega = c_1/c$. 
Hence $c \in \mathcal{C}^\prime$  determines  $\omega$.  
In view of this fact, 
we identify the moduli space  $\mathcal{M}^\prime$  with 
the region  $\mathcal{C}^\prime$.  

\begin{figure}[hbt]
  \begin{center}
    \begin{minipage}{170pt}
    \begin{picture}(140,70)
      \put(0,-7)
	{\footnotesize
        \put(0,0){\includegraphics{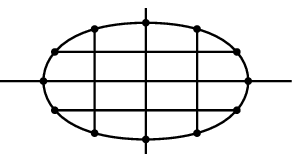}}
        \put(123,35){\makebox(0,0)[lb]{$p_1$}}
        \put(119,53){\makebox(0,0)[l]{$p_2$}}
        \put(98,63){\makebox(0,0)[lb]{$p_3$}}
        \put(77,65){\makebox(0,0)[b]{$p_4$}}
        \put(49,63){\makebox(0,0)[rb]{$p_5$}}
        \put(25,53){\makebox(0,0)[r]{$p_6$}}
        \put(20,30){\makebox(0,0)[r]{$p_7$}}
        \put(25,18){\makebox(0,0)[r]{$p_8$}}
        \put(49,7){\makebox(0,0)[rt]{$p_9$}}
        \put(70,5){\makebox(0,0)[rt]{$p_{10}$}}
        \put(91,7){\makebox(0,0)[lt]{$p_{11}$}}
        \put(119,21){\makebox(0,0)[l]{$p_{12}$}}
        \put(69,34){\makebox(0,0)[tr]{$0$}}
      }
    \end{picture}
    \caption[levelset5-1]{\\~~12 points on $L_s$}
    \label{Fig:TwelvePointsOnLevelSet}
    \end{minipage}
    \begin{minipage}{170pt}
    \begin{picture}(140,70)
      \put(0,-7)
	{\footnotesize
        \put(0,0){\includegraphics{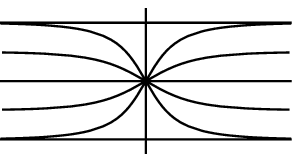}}
        \put(144,35){\makebox(0,0)[l]{$l_1$}}
        \put(144,49){\makebox(0,0)[l]{$l_2$}}
        \put(144,63){\makebox(0,0)[l]{$l_3$}}
        \put(72,69){\makebox(0,0)[l]{$l_4$}}
        \put(-4,63){\makebox(0,0)[r]{$l_5$}}
        \put(-4,49){\makebox(0,0)[r]{$l_6$}}
        \put(-4,35){\makebox(0,0)[r]{$l_7$}}
        \put(-4,21){\makebox(0,0)[r]{$l_8$}}
        \put(-4,7){\makebox(0,0)[r]{$l_9$}}
        \put(70,2){\makebox(0,0)[r]{$l_{10}$}}
        \put(144,7){\makebox(0,0)[l]{$l_{11}$}}
        \put(144,21){\makebox(0,0)[l]{$l_{12}$}}
        \put(105,42){\makebox(0,0)[]{$C_1$}}
        \put(98,56){\makebox(0,0)[t]{$C_2$}}
        \put(77,60){\makebox(0,0)[t]{$C_3$}}
        \put(63,60){\makebox(0,0)[t]{$C_4$}}
        \put(42,56){\makebox(0,0)[t]{$C_5$}}
        \put(35,42){\makebox(0,0)[]{$C_6$}}
        \put(35,28){\makebox(0,0)[]{$C_7$}}
        \put(42,14){\makebox(0,0)[b]{$C_8$}}
        \put(63,11){\makebox(0,0)[b]{$C_9$}}
        \put(79,11){\makebox(0,0)[b]{$C_{10}$}}
        \put(98,14){\makebox(0,0)[b]{$C_{11}$}}
        \put(105,28){\makebox(0,0)[]{$C_{12}$}}
      }
    \end{picture}
    \caption[levelset5-2]{\\~~12 regions on $\mathcal{C}$ }
    \label{Fig:TwelveRegions}
    \end{minipage}
  \end{center}
\end{figure}

In the level set $L_s$ for $0<s<1$,
there are twelve special points depicted in 
Figure \ref{Fig:TwelvePointsOnLevelSet} 
which give rise to pairs on  $\mathcal{M}^\prime$  rather obviously.  
Here,  
$p_1$ and $p_7$ are intersections with the real axis.  
Similarly $p_4$, $p_{10}$ lie on the imaginary axis.
$p_i$ for $i=2,3,5,6,8,9,11,12$ are points 
such that $p_3 p_{11}$, $p_5 p_9$, $p_2 p_6$ and $p_8 p_{12}$ are 
the bisectors of the segments $0 p_1$, $0 p_4$, $0 p_7$ and $0 p_{10}$ 
where  $0$  is the origin in  $\mathcal{C} \subset \mathbb{C}$.
Three points  $p_{i}$, $p_{i+2}$, $p_{i+4} ~({\rm mod}~12)$ 
together with the origin  $0$  form 
a parallelogram and give rise to a pair on  $\mathcal{M}^\prime$.  
As  $s$  decreases from  $1$  toward  $0$, 
an ellipselike closed curve  $L_s$  shrinks down to a point  $L_0$, 
however the order of $p_i$'s on the level set never changes.
Hence the trace of  $p_j$  defines a locus  $l_j$  in 
the region $\mathcal{C}$, 
and  $\mathcal{C}$  is divided by these loci into twelve regions 
denoted by $C_1$, $\cdots$, 
$C_{12}$ as in Figure \ref{Fig:TwelveRegions}.

To each pair  $(\omega,c)$  on the moduli space  $\mathcal{M}^\prime$, 
we choose four points  $0, c, c\omega, c(\omega-1)$  which 
form a parallelogram. 
Choosing the origin and the third and fourth vertices 
$c\omega, c(\omega -1)$, 
and adding the negative of the second one at the end, 
we obtain a next parallelogram spanned by  $0, c\omega, c(\omega -1), -c$. 
Doing the same shift twice, 
we obtain parallelograms spanned by  $0, c(\omega-1), -c, -c\omega$  
and  $0, -c, -c\omega, -c(\omega-1)$. 
The last one is symmetric to the original parallelogram about the origin.  
This shows that the moduli space  $\mathcal{M}^\prime = \mathcal{C}^\prime$  has 
a rotational symmetry of order  $6$.
On the other hand, 
the complex conjugation gives an orientation reversing 
involutive symmetry. 
Thus by these all together, 
the moduli space  $\mathcal{M}^\prime = \mathcal {C}^\prime$  admits 
a dihedral symmetry of order  $12$.  

Regarding  $\alpha, \beta$  as a meridian-longitude pair for  $N$, 
we discuss the hyperbolic Dehn filling coefficients 
of  $N_{\omega,c}$. 
Recall that  $\rho(\alpha)$  is a similarity defined 
by  $z \mapsto e^c z$  and 
$\rho(\beta)$  by  $z \mapsto e^{c\omega} z$ for nonzero $c$.  
Then the Dehn filling coefficients  $(\mu, \lambda) \in \mathbb{R}^2$  for 
$N_{\omega,c}$  is by definition, 
see  \cite{Thu},  
a continuous solution of the equation
\begin{equation*} 
	\rho(\alpha)^{\mu} \rho(\beta)^{\lambda} = 
		\text{rotation by $\pm 2\pi$}    
\end{equation*} 
under the normalization that  
$\rho(\alpha), \rho(\beta)$  go to the identity when 
the structure approaches the complete one.  
More precisely for our choice of  $c$  and  $\omega$, 
the Dehn filling coefficient  $(\mu, \lambda)$  of  $N_{\omega,c}$  is 
a solution of the linear equation
\begin{equation}\label{Eq:Coefficient} 
	\mu c + \lambda c \omega = \pm 2\pi i. 
\end{equation} 
Note that a pair $(\omega,c)$ corresponds to $\pm (\mu,\lambda)$ and conversely a pair $(\mu,\lambda)$ to $(\omega,\pm c)$.
We define $(\mu,\lambda)$ for $c=0$ as $(\infty,\infty)$.

To each finite Dehn filling coefficient  $(\mu, \lambda)$, 
we assign a slopelike number  $t = \mu/(\mu + \lambda)$,  
where the slope is  $\mu/\lambda$  by definition.
Though $(\omega,c)$ corresponds to coordinates $\pm (\mu,\lambda)$, the unique $t$ is assigned for them.
Dividing (\ref{Eq:Coefficient}) on both sides by  $\mu+\lambda$  and 
taking the real parts, 
we obtain an obvious identity, 
\begin{equation*} 
	{\rm Re}(t c +(1-t)c\omega) = 0,    
\end{equation*}    
which shows how to read off  $t$  with respect to the 
Dehn filling coefficient of  $N_{\omega,c}$  by 
an elementary plane geometry in  $\mathbb{C}$  as follows.  
$c$  and  $c\omega$  lie on the same level set 
of  $f$  in  $\mathcal{C}^\prime$.  
Here we draw the straight line through  $c$  and  $c\omega$.  
Then the intersection of the line with the imaginary axis 
is the point which divides the segment connecting  $c$  and  $c\omega$  
in  $1-t$  to  $t$. 

For our convenience, 
we define a function  
$T : \mathcal{C}^\prime  \to \mathbb{R} \cup \{\infty\}$  
by assigning  $t = \mu/(\mu + \lambda)$  to 
each  $c \in \mathcal{C}^\prime = \mathcal{M}^\prime$.  
It is easy to see that 
the level set of  $T$  at  $t=-1,0,1/2,1,2,\infty$  coincide with 
the union of $l_1$ and $l_7$, the union of $l_2$ and $l_8$, $\cdots$, 
the union of $l_6$ and $l_{12}$.
When  $c$  moves along the level set of  $f$  in 
counterclockwise direction,
$c\omega$  does the same.
Hence  $T$  is strictly increasing along 
a counterclockwise move on the level set of  $f$  except 
a gap at  $t = \infty$.  
\begin{figure}[hbt]
 \begin{center}
    \begin{picture}(140,70)
    \put(0,-15){
      \put(0,10){\includegraphics{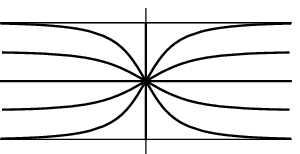}}
   \put(60,85){$t=1$}
   \put(145,72){$t=1/2$}
   \put(145,58){$t=0$}
   \put(145,44){$t=-1$}
   \put(145,30){$t=\infty$}
   \put(145,16){$t=2$}
   \put(60,0){$t=1$}
   \put(-28,72){$t=2$}
   \put(-32,58){$t=\infty$}
   \put(-34,44){$t=-1$}
   \put(-28,30){$t=0$}
   \put(-37,16){$t=1/2$}
    }
    \end{picture}
 \end{center}
  \caption[]{}
  \label{Fig:}
\end{figure}

Choose a point  $c$  in the region  $C_3$,  
then  $t$  assigned lies in  $(1/2,1)$.   
When  $c$  moves along the level set of  $T$  toward the boundary 
of $\mathcal{C}$,  
$c\omega$  must move along the level set of  $T$  in  $C_5$.
Since  ${\rm Re} (c\omega)$  goes to  $-\infty$  and   
${\rm Re}( t c + (1-t) c\omega)$  remains zero, 
${\rm Re} (c)$  must go to  $\infty$.   
The limit of the curve  $l_3$  is  $\infty+\pi i$,
and this gives ${\rm Im}(c)\rightarrow \pi$.
On the other hand, by (\ref{Eq:ConditionForAdmittingPacking}), $\cos {\rm Im}(c\omega/2)=\cos {\rm Im}(c(\omega-1)/2) \cosh {\rm Re}(c\omega/2)/\cosh {\rm Re}(c(\omega-1)/2)$.
Since $\cos {\rm Im}(c(\omega-1)/2)$ is bounded and $\cosh {\rm Re}(c\omega/2)/\cosh {\rm Re}(c(\omega-1)/2)$ is asymptotically $e^{-{\rm Re}c/2}$, $\cos {\rm Im}(c\omega/2)$ converges to $0$ and the limit of ${\rm Im}(c\omega)$ is $\pi$.
Thus,  
$t c + (1-t) c\omega \rightarrow \pi i$, 
and we have 
\begin{equation*} 
	2t c + 2(1-t) c\omega \rightarrow 2\pi i. 
\end{equation*} 
This shows that  $(\mu,\lambda)\rightarrow \pm(2t,2(1-t))$  and  
$\mu+\lambda \rightarrow \pm2$.  

\begin{figure}[hbt]
  \begin{center}
    \begin{picture}(115,110)
    \put(0,-10){
       \put(0,0){\scalebox{0.7}{\includegraphics{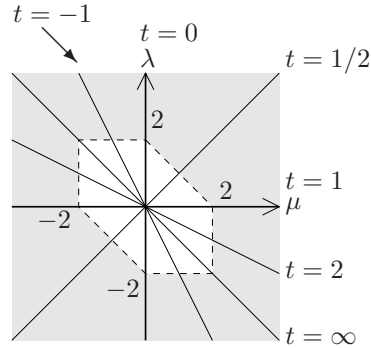}}}
       \put(2,122){$t=-1$}
       \put(13,120){\scalebox{1.3}{\vector(1,-1){10}}}
       \put(50,115){$t=0$}
       \put(105,105){$t=1/2$}
       \put(105,59){$t=1$}
       \put(105,25){$t=2$}
       \put(105,0){$t= \infty$}
       \put(105,51){$\mu$}
       \put(50,105){$\lambda$}
       \put(54,82){$2$}
       \put(80,55){$2$}
       \put(37,18){$-2$}
       \put(11,44){$-2$}
    }
    \end{picture}
    \caption[]{Dehn filling space}
    \label{Fig:Dehn Surgery Space}
  \end{center}
\end{figure}

By the dihedral symmetry of the moduli space, 
we can calculate the limit for other  $t$.
The image of the Dehn filling coefficients  $(\mu,\lambda)$  
parameterized by the moduli space  $\mathcal{M} = \mathcal{C}$  
then covers the shaded region  $\mathcal{D}$  in 
Figure \ref{Fig:Dehn Surgery Space}, 
and the boundary consists of a hexagon as follows: 
\begin{multline*}
	\{(\mu,\lambda)~|~\mu+\lambda=2, \mu \geq 0, \lambda \geq 0\} 
	\cup 	\{(\mu,\lambda)~|~\mu=2, -2 \leq \lambda \leq 0\} \\
	\cup 	\{(\mu,\lambda)~|~\lambda=-2, 0 \leq \mu \leq 2\} 
	\cup	\{(\mu,\lambda)~|~\mu+\lambda=-2, 
				\mu \leq 0, \lambda \leq 0\} \\ 
	\cup	\{(\mu,\lambda)~|~\mu=-2, 0 \leq \lambda \leq 2\}
	\cup	\{(\mu,\lambda)~|~\lambda=2, -2 \leq \mu \leq 0\}.
\end{multline*} 

%
%

\section{Degenerations}

We check how the degeneration occurs when 
they approach the boundary. 

We first think of the deformations along the path 
defined by the level set of  $T$  at  $t=1$.  
Then the deformations are realized by hyperbolic cone-manifolds and 
the meridian loop  $\alpha$  bounds a disk hitting the cone singularity.  
See  \cite{CHK}  for details about cone-manifolds.  
Since 
\begin{equation*}
	|{\rm Im}(c)|\rightarrow\pi \quad \text{and} \quad 
	|{\rm Re}(c\omega)|\rightarrow\infty,
\end{equation*}
the cone angle converges to $\pi$ and 
the length of the singular locus diverges.  
In fact, 
the Dehn filled manifold splits along 
an euclidean sub cone-manifold, 
a disk with two corner of angle  $\pi/2$  and 
a cone point of cone angle  $\pi$, 
at the limit.  
The same result holds when  $t=0$  and  $t=\infty$  
because of the dihedral symmetry.

More generally, 
choose coprime integers $p$ and $q$  and 
think of the deformations of  $N$  along the path 
on the level set of  $T$  at  $t = p/(p+q)$.  
Then again the deformations are realized by cone-manifolds and 
a loop homotopic to  $\alpha^p \beta^q$  bounds 
a disk hitting the cone singularity. 
For $t \in (0,1)$, the level set of $T$ lies in $C_2\cup l_3\cup C_3$.
When the deformations approach the boundary along this path,
${\rm Re}(c)\rightarrow\infty$, 
${\rm Im}(c)\rightarrow\pi$,
${\rm Re}(c\omega)\rightarrow -\infty$  and 
${\rm Im}(c\omega)\rightarrow \pi$.
The calculation similar to one in the previous section shows that 
\begin{equation*}
	|{\rm Im}(pc+qc\omega)| \rightarrow |(p+q)\pi|, 
\end{equation*} 
and hence the cone angle converges to $(p+q)\pi$.  
To see the length of the singular locus, 
we choose integers  $r$  and  $s$  such that $ps-qr=1$.
Then, 
$\alpha^r \beta^s$  represents a curve which intersects  $\alpha^p\beta^q$  
once and hence the translation factor of  $\rho(\alpha^r\beta^s)$, 
or equivalently  ${\rm Re}(rc+sc\omega)$  represents 
the length of the singular locus.  
Making use of the property  ${\rm Re}(pc + qc\omega) = 0$, 
we have     
\begin{equation*} 
	|{\rm Re}(rc+sc\omega)|=|\frac{rq-ps}{q}{\rm Re}(c)|\rightarrow \infty,
\end{equation*} 
which implies that the length of the singular locus diverges 
and a corresponding cone-manifold splits as before. 
We can see similar facts for other rational  $t$  
by the dihedral symmetry.  
In fact,  
the cone angle approaches  $|p|\pi$  if  $t > 1$  
and  $|q|\pi$  if  $t < 0$.  
In either case, 
the length of the singular locus diverges.

%
%

\section{Remark}

If we relax the complex affine structure modeled on  $\mathbb{C}$  to 
the complex projective structure modeled on the Riemann sphere, 
the relaxation defines the involution  $\tau$  
on  $\mathbb{H}^2 \times \mathbb{C}$  
which sends  $(\omega, c)$  to  $(\omega, -c)$.  
The moduli space  $\mathcal{M}_*$  of 
complex projective tori 
admitting a one-circle packing is 
the quotient of  $\mathcal{M}$  by  $\tau$.   
This space is studied in \cite{KMT}  by the cross ratio parameters  
and represented as a semi algebraic set by  
\begin{equation*} 
	\mathcal{M}_* =\{(x,y,z) \, | \, xyz-x-y-z=0,x,y,z>0\}, 
\end{equation*}
see the appendix in  \cite{KMT}.  

Assigning the holonomy representation to each 
projective structure, 
we obtain a map 
\begin{equation*} 
	\mathcal{M}_* \rightarrow \mathbb{H}^2 \times \mathbb{C}/\tau
		\rightarrow  X(\pi_1), 
\end{equation*}
of  $\mathcal{M}_*$  to 
the space  $X(\pi_1)$   of representations
of the fundamental group of the torus in  ${\rm PGL}(2,\mathbb{C})$  
up to conjugation.  
It is a regular map with respect to the algebraic 
structures on the source and target, 
and easily seen to be a proper injective homeomorphism 
by concrete computation in  \cite{KMT}.  
This in particular implies also that 
the boundary we described here is the true one 
since holonomy representations diverge at the boundary.

\end{document}